\documentclass[11pt,letterpaper,reqno]{amsart}
\usepackage{tikz}
\usetikzlibrary{positioning, shapes.geometric, arrows.meta}
\usepackage{amssymb}
\usepackage{amsmath}
\usepackage{amsthm}
\usepackage{amsfonts}
\usepackage{bbm}
\usepackage{graphicx}
\usepackage[T1]{fontenc}
\usepackage{doi}
\usepackage{bookmark}
\usepackage{hyperref}
\hypersetup{pdfstartview={FitH}}

\newtheorem{thm}{Theorem}[section]
\newtheorem{lem}[thm]{Lemma}

\newtheorem{cor}[thm]{Corollary}

\newtheorem{ques}[thm]{Question}
\newtheorem{prob}[thm]{Problem}
\newtheorem{conj}[thm]{Conjecture}
\theoremstyle{definition}
\newtheorem{remark}[thm]{Remark}
\newtheorem{defn}[thm]{Definition}

\numberwithin{equation}{section}

\makeatother

\begin{document}

\title[On the Positive and Negative $p$-Energies of Graphs under Edge Addition]{On the Positive and Negative $p$-Energies of Graphs under Edge Addition}

\author[Quanyu Tang]{Quanyu Tang}
\author[Yinchen Liu]{Yinchen Liu}
\author[Wei Wang]{Wei Wang}

\address{School of Mathematics and Statistics, Xi'an Jiaotong University, Xi'an 710049, P. R. China}
\email{tang\_quanyu@163.com}
\address{Institute for Interdisciplinary Information Sciences, Tsinghua University, Beijing 100084, P. R. China}
\email{liuyinch23@mails.tsinghua.edu.cn}
\address{School of Mathematics and Statistics, Xi'an Jiaotong University, Xi'an 710049, P. R. China}
\email{wang\_weiw@xjtu.edu.cn}


\begin{abstract} 
In this paper, we introduce the concepts of positive and negative $p$-energies
of graphs and investigate their behavior under edge addition. Specifically, we generalize the classical notions of positive and negative square energies to the $p$-energy setting, denoted by $\mathcal{E}_p^{+}(G)$ and $\mathcal{E}_p^{-}(G)$, respectively. We establish improved lower bounds for these quantities under edge addition, which sharpen existing results by Abiad et al.\ in the case $p=2$. Furthermore, we address the monotonicity problem for $\mathcal{E}_p^{+}(G)$ under edge addition, and construct a family of counterexamples showing that monotonicity fails for $1 \leq p < 3$. Finally, we conclude with several open problems for further investigation.
\end{abstract}

\maketitle

\section{Introduction}

Let us begin with some definitions and notation. Throughout this paper, we consider only graphs that are simple (i.e., without loops or multiple edges), undirected, and unweighted. Let \( G = (V, E) \) be a graph of order \( n \) and size \( m \). The \emph{adjacency matrix} of \( G \) is the \( n \times n \) matrix \( A(G) = [a_{ij}] \), where \( a_{ij} = 1 \) if vertices \( v_i \) and \( v_j \) are adjacent, and \( a_{ij} = 0 \) otherwise. The \emph{eigenvalues} of \( G \) are the eigenvalues of \( A(G) \). Since \( A(G) \) is a real symmetric matrix, all of its eigenvalues are real.

Let $s^{+}(G)$ ($s^{-}(G)$) be the sum of the squares of the positive (negative) eigenvalues of the adjacency matrix $A(G)$ of $G$. We refer to $s^{+}(G)$ and $s^{-}(G)$ as the \emph{positive and negative square energies} of the graph. In \cite{EFGW16}, Elphick, Farber, Goldberg, and Wocjan  proposed the following conjecture:

\begin{conj}[Elphick et al.~\cite{EFGW16}] \label{conj:extension-Hong}
Let \( G \) be a connected graph with \( n \) vertices. Then
\[
\min\{s^+(G), s^-(G)\} \geq n - 1.
\]
\end{conj}

The above conjecture has attracted some attention from the spectral graph theory community; for example, it is listed as the first conjecture in a recent survey \cite{LiuNing2023} by Liu and Ning. Conjecture~\ref{conj:extension-Hong} has been verified for several graph classes, including bipartite graphs, regular graphs, complete $q$-partite graphs, hyper-energetic graphs, and barbell graphs (see \cite{Abiad2023, ElphickLinz2024} for other partial results). Recently, Zhang \cite{Zhang2024} and Akbari et al. \cite{Akbari2024} have made substantial progress on this conjecture, both utilizing the super-additivity of the square energies. They proved, respectively, that \( \min\{s^+(G), s^-(G)\} \geq n - \gamma \) and \( \min\{s^+(G), s^-(G)\} \geq \frac{3n}{4} \), where \(\gamma\) is the domination number of \(G\). However, Conjecture \ref{conj:extension-Hong} still appears to be highly challenging.

Since then, the sum of the squares of the positive and negative eigenvalues has been a topic of interest. In particular, the sum of the squares of the positive eigenvalues, $s^+(G)$, seems to exhibit more favorable properties than that of the negative eigenvalues, $s^-(G)$. Therefore, in Cioaba et al. \cite{Cioaba}, Guo posed the following conjecture during the workshop open problem session, which generalizes to $s^+$ a known result for the spectral radius, namely that $\rho(G + uv) \geq \rho(G)$:

\begin{conj}[\cite{Cioaba}]\label{conj2}
If $G$ is a graph and $uv \notin E(G)$, then
\[
s^+(G + uv) \geq s^+(G).
\]
\end{conj}

Unfortunately, Conjecture \ref{conj2} turns out to be false; see e.g. Abiad et al. \cite[Example 2.6]{Abiad2023}. Additionally, Abiad et al. \cite{Abiad2023} also established lower bounds for \( s^{\pm}(G) \) when an edge is removed from the graph:

\begin{thm}[\cite{Abiad2023}, Theorem 2.5]\label{Abiad}
 Let $G$ be a graph, let $H = G - e$ (where $e$ is an edge of $G$), and let $\lambda_1 \geq \cdots \geq \lambda_n $ and $\theta_1 \geq \cdots \geq \theta_n$ denote the eigenvalues of $A(G)$ and $A(H)$, respectively. If $H$ has at least two positive eigenvalues and at least two negative eigenvalues, then $s^+(G) \geq s^+(H) - \theta_{2}^2$ and $s^-(G) \geq s^-(H) - \theta_{n}^2$.
\end{thm}

In Nikiforov \cite{Nikiforov2012}, the Schatten $p$-norms of graphs were studied. We now recall the definition of the Schatten norms for square matrices:

\begin{defn}

Let \(A\) be an \(n \times n\) complex matrix, and let \(p \geq 1\). The \emph{Schatten \(p\)-norm} \(\|A\|_{p}\) is defined as \[
\|A\|_{p} := \left( \sum_{i=1}^n \sigma_i(A)^p \right)^{1/p},
\] where \(\sigma_1(A), \dots, \sigma_n(A)\) are the singular values of \(A\).
\end{defn}

If $G$ is a graph with adjacency matrix $A(G)$, for short we write $\|G\|_{p}$ for $\|A(G)\|_{p}$. Let \( \lambda_1, \dots, \lambda_n \) be the eigenvalues of \( A(G) \), collectively known as the \emph{spectrum} of \( G \). Then the \textit{$p$-energy} of $G$ is defined as the sum of the \( p \)th power of the absolute values of the eigenvalues of its adjacency matrix$$\mathcal{E}_p(G)=\sum^n_{i=1}|\lambda_i|^p,$$where $p$ is a positive real number. Clearly, the $p$-energy $\mathcal{E}_p(G)$ of a graph $G$ is exactly $\|G\|_{p}^p$. When \( p = 1 \), \( \mathcal{E}_1(G) \) is simply denoted as \( \mathcal{E}(G) \), known as the \textit{energy} of \( G \). Originating from applications in molecular chemistry, graph energy has attracted considerable attention. The \( p \)-energy is also known as the \textit{\( p \)-Schatten energy} in other works, such as Arizmendi and Guerrero \cite{Arizmendi23}. In this paper, we follow the terminology used by Akbari et al. \cite{Akbari2020}.

Inspired by the aforementioned notions, in this paper we generalize these concepts and investigate the sums of the absolute values of the positive and negative eigenvalues raised to the $p$th power. It is worth noting that some tools from matrix analysis can be readily applied to derive general results for $p \geq 1$. To formalize this, let $\mathcal{E}_p^+(G)$ ($\mathcal{E}_p^-(G)$) denote the sum of the absolute values of the positive (negative) eigenvalues raised to the $p$th power of the adjacency matrix $A(G)$ of a graph $G$, i.e.

\[
\mathcal{E}_p^+(G) = \sum_{\lambda_i > 0} \lambda_i(G)^p \quad \text{and} \quad \mathcal{E}_p^-(G) = \sum_{\lambda_i < 0} |\lambda_i(G)|^p.
\] The parameters $\mathcal{E}_p^+(G)$ and $\mathcal{E}_p^-(G)$ are called the \emph{positive $p$-energy} and the \emph{negative $p$-energy} of $G$, respectively. Clearly, $\mathcal{E}_p^+(G) + \mathcal{E}_p^-(G) = \mathcal{E}_p(G)$. In particular, when $p = 2$, we write $\mathcal{E}_2^+(G)$ and $\mathcal{E}_2^-(G)$ as $s^+(G)$ and $s^-(G)$, respectively, and refer to them as the \emph{positive square energy} and \emph{negative square energy} of $G$. When $p = 1$, it is well known that $\mathcal{E}_1^+(G) = \mathcal{E}_1^-(G) = \frac{1}{2} \mathcal{E}(G)$.

Naturally, we may ask whether a similar edge-addition monotonicity holds for the general $\mathcal{E}_p^+(G)$, as in Conjecture \ref{conj2}:

\begin{ques}\label{q1}
Given $p \geq 1$, let $G$ be a graph and suppose that $uv \notin E(G)$. Does it always hold that $\mathcal{E}_p^+(G + uv) \geq \mathcal{E}_p^+(G)$?
\end{ques}

\begin{remark}
After this paper was made publicly available but before its formal journal publication, some progress had already been made on the conjecture proposed in its final section. For instance, see the recent work of Akbari, Kumar, Mohar, and Pragada~\cite{Akbari2025}. 
\end{remark}
\begin{remark}
We also remark that the study of positive and negative \(p\)-energies is not merely a theoretical generalization. In a recent work~\cite{Tang2025}, the first author established novel lower bounds on the chromatic number \(\chi(G)\) using positive and negative \(p\)-energies, which not only generalize but also strengthen the classical spectral bounds derived from graph energy, such as those by Ando and Lin~\cite{Ando2015}, and Elphick and Wocjan~\cite{Elphick2017}. Notably, it was shown that for certain graphs, non-integer values of \(p\) provide sharper lower bounds on \(\chi(G)\) than all previously known spectral bounds. This demonstrates that the parameters \(\mathcal{E}_p^+(G)\) and \(\mathcal{E}_p^-(G)\) have practical value in spectral graph theory beyond the classical \(p = 1\) or \(2\) settings.
\end{remark}

The paper is organized as follows. In Section~\ref{sec2}, we investigate the behavior of $\mathcal{E}_p^{\pm}(G)$ under edge addition. As a corollary, we present tighter bounds than those in Abiad et al.~\cite[Theorem~2.5]{Abiad2023}, improving the order of the difference in positive and negative square energies before and after edge addition. In Section~\ref{sec3}, we address part of Question~\ref{q1}, proving that $\mathcal{E}_p^+(G)$ is not monotonically increasing under edge addition for $1 \leq p < 3$. In doing so, we construct a family of graphs that likely represent the most extremal cases, offering a promising direction for further investigation of Question~\ref{q1}. Finally, in Section~\ref{sec4}, we provide some open problems for future study.

\section{Improved Bounds under Edge Addition}\label{sec2}

In this section, we analyze the effect on the positive and negative $p$-energies $\mathcal{E}_p^{\pm}(G)$ when adding an edge to a graph $G$. Note that Conjecture \ref{conj:extension-Hong} holds for graphs with exactly one positive eigenvalue or exactly
one negative eigenvalue \cite{EFGW16}, hence the assumption that there are at least two positive and two negative eigenvalues does not diminish the relevance of the result.

We present an analogue of Theorem \ref{Abiad} about $\mathcal{E}_p^{\pm}(G)$, which is tighter compared to Theorem \ref{Abiad}. For easier comparison, we express the results in this section in terms of removing an edge.

\begin{thm}\label{t2}
Let $p \geq 1$. Let $G$ be a graph, let $H = G - e$ (where $e$ is an edge of $G$), and let $\lambda_1 \geq \cdots \geq \lambda_n$ and $\theta_1 \geq \cdots \geq \theta_n$ denote the eigenvalues of $A(G)$ and $A(H)$, respectively. If $H$ has at least two positive
eigenvalues and at least two negative eigenvalues, then we have the following inequalities:
\[
\mathcal{E}_p^+(G) \geq \mathcal{E}_p^+(H) + \max\{\theta_2 - 1, 0\}^p - \theta_2^p,
\]and\[
\mathcal{E}_p^-(G) \geq \mathcal{E}_p^-(H) + \max\{-\theta_{n}-1, 0\}^p - |\theta_n|^p.
\]
\end{thm}

Set $p=2$ in Theorem \ref{t2}, then we can directly obtain the following:

\begin{cor}\label{c1}
Let $G$ be a graph, let $H = G - e$ (where $e$ is an edge of $G$), and let $\lambda_1 \geq \cdots \geq \lambda_n $ and $\theta_1\geq \cdots \geq \theta_{n}$ denote the eigenvalues of $A(G)$ and $A(H)$, respectively. If $H$ has at least two positive
eigenvalues and at least two negative eigenvalues, then we have the following inequalities:\[
s^+(G) \geq
\begin{cases}
s^+(H) - \theta_2^2,&  \text{if } \theta_2 < 1 \\
s^+(H) - 2\theta_2 + 1,
&  \text{if } \theta_2 \geq 1
\end{cases},
s^-(G) \geq
\begin{cases}
s^-(H) - \theta_n^2,& \text{if } \theta_n > -1 \\
s^-(H) + 2\theta_n + 1,&  \text{if } \theta_n \leq -1
\end{cases}.
\]
\end{cor}

\begin{remark} Clearly, if we consider lower bounds for $s^+(G)$, it follows from $\theta_2^2- 2\theta_2 + 1=(\theta_2-1)^2 \geq 0$ that Corollary \ref{c1} is tighter compared with Theorem \ref{Abiad}. Similarity, if we consider lower bounds for $s^{-}(G)$, it follows from $\theta_n^2+ 2\theta_n + 1=(\theta_n+1)^2 \geq 0$ that Corollary \ref{c1} is tighter as well. The lower bound in Theorem \ref{Abiad} is of second order, while our result provides a first-order lower bound, which evidently improves the order.
\end{remark}

As a partial progress toward Conjecture~\ref{conj:extension-Hong}, we obtain the following result as another immediate corollary.

\begin{cor}
If \( G \) is a graph on \( n \) vertices and \( H = G - e \) (where \( e \) is an edge of \( G \)) satisfies 
\(
s^{+}(H) + \max\{\theta_2 - 1, 0\}^2 - \theta_2^2 \geq n - 1 \) and \(s^{-}(H) + \max\{-\theta_n - 1, 0\}^2 - |\theta_n|^2 \geq n - 1,
\)
then 
\(
\min\{s^{+}(G), s^{-}(G)\} \geq n - 1.
\)
\end{cor}

The proof of our result, like that of Abiad et al. \cite{Abiad2023}, also relies on a result on edge interlacing for the adjacency matrix established by Hall, Patel, and Stewart
 in \cite[Theorem 3.9]{Hall2009}:

\begin{lem}[\cite{Hall2009}]\label{Hall} Let $G$ be a simple graph and let $H = G - e$, where $e$ is an edge of $G$. Let $\lambda_1 \geq \lambda_2 \geq \cdots \geq \lambda_n$ and $\theta_1 \geq \theta_2 \geq \cdots \geq \theta_n$ be the eigenvalues of the adjacency matrices $A(G)$ and $A(H)$, respectively. Then, for each $i = 2, 3, \dots, n-1$, we have
\[
\lambda_{i-1} \geq \theta_i \geq \lambda_{i+1},
\]
with $\theta_1 \geq \lambda_2$ and $\theta_n \leq \lambda_{n-1}$.

\end{lem}

Another interesting method of comparing the distribution of two sets of real numbers is the concept of \emph{majorization}. It may be introduced in the following way
 (see e.g.~\cite[Ch.~1, Sect.~A]{MarshallOlkin1979}). For any vector $\mathbf{x} = (x_1, \dots, x_n) \in \mathbb{R}^n$, we denote by $\mathbf{x}^\downarrow=(x_{1}^\downarrow, \dots, x_{n}^\downarrow)$ a rearrangement of the components of $\mathbf{x}$ such that
\[
x_{1}^\downarrow \geq x_{2}^\downarrow \geq \dots \geq x_{n}^\downarrow.
\]

\begin{defn}
For any two vectors $\mathbf{x} = (x_1, \dots, x_n)$ and $\mathbf{y} = (y_1, \dots, y_n)$ from $\mathbb{R}^n$, we say that $\mathbf{y}$ \emph{weakly majorizes} $\mathbf{x}$, and write this as $\mathbf{y} \succ_w \mathbf{x}$, if
\[
\sum_{k=1}^{t} y_{k}^\downarrow \geq \sum_{k=1}^{t} x_{k}^\downarrow, \quad t = 1, \dots, n.
\]
Furthermore, we say that $\mathbf{y}$ \emph{(strongly) majorizes} $\mathbf{x}$, and write this as $\mathbf{y} \succ \mathbf{x}$, if, in addition, equality occurs in (3) when $t = n$.
\end{defn}

 We know that the Lidskii inequalities \cite[Theorem 4.3.47 (b)]{Horn13} are the basis for many important perturbation bounds. Therefore, in addition to Lemma \ref{Hall}, we also use Lidskii inequalities, which have not been previously utilized in the analysis of $s^+(G)$ under edge addition. Here, we employ an equivalent version of the Lidskii inequalities, which can be found in \cite[4.3. P24]{Horn13}:

\begin{lem}\label{Lidskii}
Let $A, B \in M_n$ be Hermitian. Let $\lambda(A)$, $\lambda(B)$, and $\lambda(A - B)$, respectively, denote the real $n$-vectors of eigenvalues of $A$, $B$, and $A - B$. Then
\[
    \lambda(A - B) ~\text{majorizes}~ \lambda(A)^\downarrow - \lambda(B)^\downarrow.
\]
\end{lem}

We are now ready to present the proof of Theorem \ref{t2}.

\begin{proof}[Proof of Theorem \ref{t2}]
We assume that $H$ has $r$ positive eigenvalues, where $r \geq 2$. Now we claim that
\[
\left( \lambda_1, \max\{\lambda_2, 0\}, \dots, \max\{\lambda_r, 0\}\right) \succ_w \left( \theta_1, \max\{\theta_2 - 1, 0\}, \theta_3, \theta_4, \dots, \theta_r\right).   \tag{1}\label{eq1}
\]
By the Perron-Frobenius theorem, it is clear that $\lambda_1 \geq \theta_1$, since $H \subseteq G$.

If $\theta_2 \leq 1$, then $\max\{\theta_2 - 1, 0\} = 0$, thus we only need to show that
\[
\lambda_1 + \sum_{k=2}^{t} \max\{\lambda_k, 0\} \geq \theta_1 + \sum_{k=3}^{t} \theta_k \tag{2}\label{eq2}
\]
holds for $2 \leq t \leq r$. Since Lemma \ref{Hall} tells us that $\lambda_{k-1} \geq \theta_k \geq 0$ for $2 \leq k \leq r$, we have
\[
\sum_{k=2}^{t} \max\{\lambda_k, 0\} = \sum_{k=3}^{t+1} \max\{\lambda_{k-1}, 0\} \geq \sum_{k=3}^{t} \lambda_{k-1} \geq \sum_{k=3}^{t} \theta_k .
\]
Hence, inequality (\ref{eq2}) holds since $\lambda_1 \geq \theta_1$.

If $\theta_2 > 1$, then $\max\{\theta_2 - 1, 0\} = \theta_2 - 1$, thus we only need to show
\[
\lambda_1 + \sum_{k=2}^{t} \max\{\lambda_k, 0\} \geq \sum_{k=1}^{t} \theta_k - 1 \tag{3}\label{eq3}
\]
holds for $2 \leq t \leq r$. Let $F$ denote the adjacency matrix of the graph obtained from the empty graph by adding an edge $e$, then $F = A(G) - A(H)$. Let $\lambda_j(X)$ denote the eigenvalues of $X$, where $1 \leq j \leq n$. Thus by Lemma \ref{Lidskii}, $\lambda(-F) \succ \lambda(A(H))^\downarrow - \lambda(A(G))^\downarrow$, which means
\[
(1, 0, \dots, 0, -1) \succ (\theta_1 - \lambda_1, \theta_2 - \lambda_2, \dots, \theta_n - \lambda_n),
\]
i.e. for $1 \leq t \leq n$, we have
\[
\sum_{j=1}^{t} (\theta_j - \lambda_j) \leq \sum_{j=1}^{t} \lambda_j(-F)^\downarrow.
\]
Hence
\[
\sum_{j=1}^{t} \lambda_j \geq \sum_{j=1}^{t} \theta_j - \sum_{j=1}^{t} \lambda_j(-F)^\downarrow.
\]
Therefore, for $2 \leq t \leq r$, we have
\[
\lambda_1 + \sum_{k=2}^{t} \max\{\lambda_k, 0\} \geq \sum_{j=1}^{t} \lambda_j \geq \sum_{j=1}^{t} \theta_j - 1,
\]
which means (\ref{eq3}) holds and thus we prove the claim of (\ref{eq1}).

By the Hardy-Littlewood-P\'{o}lya Theorem, it follows that $\mathbf{y} \succ_w \mathbf{x}$ implies $f(\mathbf{y}) \succ_w f(\mathbf{x})$ for any continuous, convex, and nondecreasing function $f$. For $p \geq 1$, we have
\[
\lambda_1^p + \sum_{k=2}^{r} \max\{\lambda_k, 0\}^p \geq \sum_{k=1}^{r} \theta_k^p + \max\{\theta_2 - 1, 0\}^p - \theta_2^p,
\]
as the majorization relation in (\ref{eq1}) holds. Hence, we have that
\[
\mathcal{E}_p^+(G) \geq \sum_{k=1}^{r} \max\{\lambda_k, 0\}^p \geq \mathcal{E}_p^+(H) + \max\{\theta_2 - 1, 0\}^p - \theta_2^p,
\]and the first inequality in Theorem \ref{t2} follows.

Next, we assume that \( H \) has \( q \) negative eigenvalues, where \( q \geq 2 \). Analogously, we obtain the following majorization relation:
\begin{align}
&(-\lambda_n, \max\{-\lambda_{n-1}, 0\}, \max\{-\lambda_{n-2}, 0\}, \dots, \max\{-\lambda_{n+1-q}, 0\}) \notag \\
&\qquad \succ_w (\max\{-\theta_n - 1, 0\}, -\theta_{n-1}, -\theta_{n-2}, \dots, -\theta_{n+1-q}). \notag
\end{align}
Hence, we have that
\[
 \sum_{j=1}^{q} \max\{-\lambda_{n+1-j}, 0\}^p \geq \max\{-\theta_{n}-1, 0\}^p + \sum_{j=2}^{q} |\theta_{n+1-j}|^p.
\]Therefore,
\[
\mathcal{E}_p^-(G) \geq \sum_{j=1}^{q} \max\{-\lambda_{n+1-j}, 0\}^p \geq \mathcal{E}_p^-(H) + \max\{-\theta_{n}-1, 0\}^p - |\theta_n|^p,
\]as desired.
\end{proof}




\section{Monotonicity for Edge Addition}\label{sec3}

In this section, we shall construct a family of graphs which show that the monotonicity property fails for edge addition when 
$1\leq p<3$; in particular, this provide a family of counterexamples to Conjecture \ref{conj2}.
Before doing so, we first recall some results from the theory of graph spectra.

\begin{defn}
Let $G$ be a graph with adjacency matrix $A(G)$. A partition $\pi=(V_1,V_2,\dots,V_k)$ of $V(G)$ is an \emph{equitable partition} if every vertex in $V_i$ has $b_{ij}$ (a constant only dependent on $i$ and $j$) neighbours in $V_j$ for all $i,j\in\{1,2,\dots,k\}$. The $k\times k$ matrix $B_\pi(G)=(b_{ij})$ is called the \emph{divisor matrix} of the partition $\pi$.
\end{defn}

The following lemma shows that the spectrum of $B_\pi(G)$ is a subset of that of $A(G)$.

\begin{lem}[\cite{Cvetkovic2010}, Theorem 3.9.5]\label{lem-equitable partition}
The characteristic polynomial of any divisor matrix $B_\pi$ of a graph $G$
divides the characteristic polynomial of $G$.
\end{lem}

In addition, we will also need the vertex version of the interlacing property of the eigenvalues of $A(G)$ and the eigenvalues of $A(G - v)$, where $v$ is a vertex of $G$.

\begin{lem}[\cite{Hall2009}, Theorem 2.3]\label{lem4}
Let $G$ be a graph and $H = G - v$, where $v$ is a vertex of $G$. If \(
\lambda_1 \geq \lambda_2 \geq \cdots \geq \lambda_n \) and \( \theta_1 \geq \theta_2 \geq \cdots \geq \theta_{n-1}
\) are the eigenvalues of $A(G)$ and $A(H)$, respectively, then
\[
\lambda_i \geq \theta_i \geq \lambda_{i+1} \quad \text{for each} \quad i = 1, 2,\dots, n-1.
\tag{6}\label{eq6}\]
\end{lem}

\begin{remark} When the inequalities (\ref{eq6}) are satisfied, we say that the eigenvalues $\theta_i$
\emph{interlace} the eigenvalues $\lambda_j$ .
\end{remark}

\begin{figure}[!h]
  \centering
\begin{tikzpicture}[x=1pt,y=-1pt,line width=5.0]
\draw[line width=1.0](66.5,37.3) ellipse (66.5 and 37.38);
\node at (66.5,36.1){$K_{n_1-1}$};
\node at (276.5,36.1){$K_{n_2-1}$};
\draw[line width=1.0](276.5,37.3) ellipse (66.5 and 37.38);
\draw(132.9,38.5)--(210,38.5);
\draw(62.7,128)--(62.7,74.6);
\draw(272.7,128)--(272.7,74.6);
\node at (63.2,148.1){$v_1$};
\node at (273.2,148.1){$v_2$};
\draw[red, fill=red](62.7,128) circle (5);
\draw[red, fill=red](272.7,128) circle (5);
\draw[fill=black](132.9,38.5) circle (2);
\draw[fill=black](210,38.5) circle (2);
\draw[fill=black](62.7,74.6) circle (2);
\draw[fill=black](272.7,74.6) circle (2);
\end{tikzpicture}
  \caption{The complement of double star graph $\overline{S_{n_1,n_2}}$.}\label{fig1}
\end{figure}
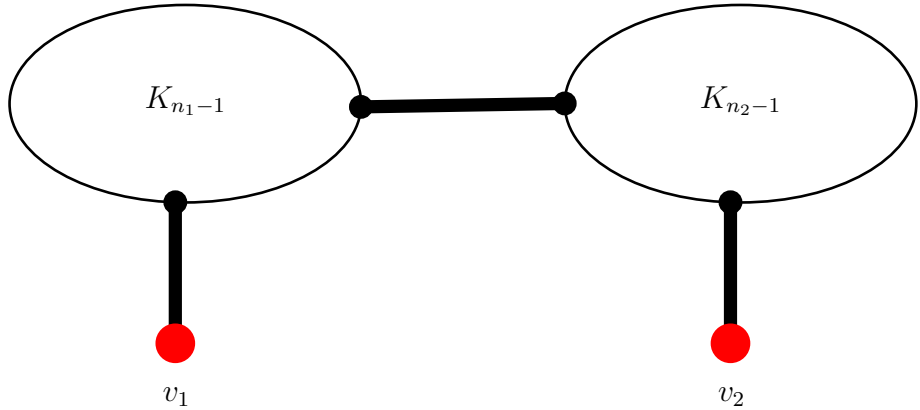

Now, we describe our construction as follows. Let $S_{n_1,n_2}$ denote the \emph{double star} graph obtained by adding an edge between the centers of the stars $K_{1,n_1-1}$ and $K_{1,n_2-1}$, where $v_1$ and $v_2$ denote the respective centers of $K_{1,n_1-1}$ and $K_{1,n_2-1}$. Let $\overline{S_{n_1,n_2}}$ be the complement of $S_{n_1,n_2}$, as illustrated in Figure~\ref{fig1}, where the bold lines indicate that all vertices in one set are adjacent to all vertices in the other.

The next theorem gives the spectrum of $\overline{S_{n,n}}$ by the technique of equitable partition and eigenvalue interlacing.

\begin{thm}\label{thm-CSmm}
The eigenvalues of graph $\overline{S_{n,n}}$ are 
\(\frac{1}{2} \left(-1+\sqrt{4n-3}\right)\), \(\frac{1}{2}\left(-1-\sqrt{4n-3}\right)\), \\\(\frac{1}{2} \left(2n-3+\sqrt{4n^2-8n+5}\right)\),
\(\frac{1}{2} \left(2n-3-\sqrt{4n^2-8n+5}\right)\), \(\underbrace{-1,\dots,-1}_{2n-4}\).
\end{thm}
\begin{proof}
Let $\pi=V_1\cup V_2\cup V_3\cup V_4$ be a partition of $V(\overline{S_{n,n}})$ given by $V_1=\{v_1\}$, $V_2=\{v_2\}$, $V_3=V(K_{n-1})$ and  $V_4=V(K_{n-1})$. Then $\pi$ is an equitable partition with divisor matrix 
\begin{equation*}
  B_\pi=\left(
\begin{array}{cccc}
 0 & 0 & 0 & n-1 \\
 0 & 0 & n-1 & 0 \\
 0 & 1 & n-2 & n-1 \\
 1 & 0 & n-1 & n-2 \\
\end{array}
\right).
\end{equation*}
The characteristic polynomial of $B_\pi$ is
\begin{equation*}
P_{B_\pi}(x)=\left(x^2-(2n-3)x-(n-1)\right)\left(x^2+x-(n-1)\right).
\end{equation*}
Therefore, by Lemma \ref{lem-equitable partition}, we know that the eigenvalues of $A(G)$ include all the roots of $P_{B_\pi}(x)$, which are \(\frac{1}{2} \left(-1+\sqrt{4n-3}\right)\), \(\frac{1}{2}\left(-1-\sqrt{4n-3}\right)\), \(\frac{1}{2} \left(2n-3+\sqrt{4n^2-8n+5}\right)\) and \(\frac{1}{2} \left(2n-3-\sqrt{4n^2-8n+5}\right)\). Since the complete graph $K_{2n-2}$ has eigenvalues $-1$ with multiplicity $2n-3$, which interlace the eigenvalues of graph $\overline{S_{n,n}}-v_1$, we see that $\overline{S_{n,n}}-v_1$ has eigenvalues $-1$ with multiplicity at least $2n-4$ by Lemma \ref{lem4}. Consider again the eigenvalue interlacing of $\overline{S_{n,n}}-v_1$ and $\overline{S_{n,n}}$, we conclude that $\overline{S_{n,n}}$ has eigenvalues $-1$ with multiplicity at least $2n-5$.
Since the sums of eigenvalues $B_\pi$ and $A(\overline{S_{n,n}})$ are $2n-4$ and $0$, respectively, we conclude that $A(\overline{S_{n,n}})$ has $-1$ as an eigenvalue with multiplicity $2n-4$.
\end{proof}

Let $(\overline{S_{n_1,n_2}})^+$ be the graph obtained from $\overline{S_{n_1,n_2}}$ by adding an edge $v_1v_2$, where $v_1$ and $v_2$ are the respective centers of $K_{1,n_1-1}$ and $K_{1,n_2-1}$ as shown in Figure~\ref{fig1}. By similar arguments as Theorem \ref{thm-CSmm}, we can derive the spectrum of $(\overline{S_{n,n}})^+$ in the following theorem.

\begin{thm}\label{thm-CSmm2}
The eigenvalues of graph $(\overline{S_{n,n}})^+$ are $n-1+\sqrt{n^2-3 n+3}$, $-1+\sqrt{n-1}$, $-1-\sqrt{n-1}$,$n-1-\sqrt{n^2-3 n+3}$, $\underbrace{-1,\dots,-1}_{2n-4}$.
\end{thm}

Now we are in a position to present the following
\begin{thm}\label{t5}
For each fixed $1 \leq p < 3$, there exists a family of connected graphs $G$ and $uv \notin E(G)$ such that
\[
\mathcal{E}_p^+(G + uv) < \mathcal{E}_p^+(G).
\]
\end{thm}

\begin{proof}
Let \( n \geq 3 \). Let \( G = \overline{S_{n,n}} \), and let \( G + uv = (\overline{S_{n,n}})^+ \) denote the graph obtained by adding an edge \( uv \) as described above. By Theorem \ref{thm-CSmm} and Theorem \ref{thm-CSmm2}, we know that the positive eigenvalues of \( \overline{S_{n,n}} \) are
\[
\lambda_1 = \frac{1}{2} \left(-1 + \sqrt{4n - 3} \right), \quad \lambda_2 = \frac{1}{2} \left( 2n - 3 + \sqrt{4n^2 - 8n + 5} \right),
\]
and the positive eigenvalues of \( (\overline{S_{n,n}})^+ \) are
\[
\theta_1 = -1 + \sqrt{n - 1}, \quad \theta_2 = n - 1 + \sqrt{n^2 - 3n + 3}, \quad \theta_3 = n - 1 - \sqrt{n^2 - 3n + 3}.
\] Let \( f(n) = \mathcal{E}_p^{+}\left( \overline{S_{n,n}} \right) - \mathcal{E}_p^{+}\left( (\overline{S_{n,n}})^+ \right) = \lambda_1^p + \lambda_2^p - \theta_1^p - \theta_2^p - \theta_3^p \). Then it suffices to prove that for any \( 1 \leq p < 3 \), \( f(n) > 0 \) holds when \( n \) is sufficiently large.

For \( p = 1 \), the result is straightforward, since
\begin{eqnarray*}
f(n) &=& \lambda_1 + \lambda_2 - \theta_1 - \theta_2 - \theta_3 \\
&=& \sqrt{(n - 1)^2 + \frac{1}{4}} - (n - 1)
+ \sqrt{n - \frac{3}{4}} - \sqrt{n - 1} > 0
\end{eqnarray*}for all \( n \geq 3 \). Therefore, we now only need to consider the case \( 1 < p < 3 \).

By the Taylor series (using Big O notation), for \( \alpha > 0 \), we have

\[
(1 + x)^\alpha = 1 + \alpha x + \frac{\alpha(\alpha - 1)}{2} x^2 + O(x^3),
\]
as \( x \to 0 \). Hence,
\[
\sqrt{n - 1} = \sqrt{n} \cdot \sqrt{1 - \frac{1}{n}} = \sqrt{n} \left( 1 + O\left(\frac{1}{n}\right) \right) = \sqrt{n} + O\left(\frac{1}{\sqrt{n}}\right).
\]
Similarly,
\[
\sqrt{4n - 3} = 2\sqrt{n} + O\left(\frac{1}{\sqrt{n}}\right),
\]
\[
\sqrt{4n^2 - 8n + 5} = 2n - 2 + \frac{1}{4n} + O\left(\frac{1}{n^3}\right),
\]
\[
\sqrt{n^2 - 3n + 3} = n - \frac{3}{2} + \frac{3}{8n} + O\left(\frac{1}{n^3}\right).
\] Thus, \begin{eqnarray*}
f(n) &=& \left( \sqrt{n} - \frac{1}{2} + O\left(\frac{1}{\sqrt{n}}\right) \right)^p + \left( 2n - \frac{5}{2} + \frac{1}{8n} + O\left(\frac{1}{n^3}\right) \right)^p
\\&&- \left( \frac{1}{2} - \frac{3}{8n} + O\left(\frac{1}{n^3}\right) \right)^p - \left( \sqrt{n} - 1 + O\left(\frac{1}{\sqrt{n}}\right) \right)^p \\&&- \left( 2n - \frac{5}{2} + \frac{3}{8n} + O\left(\frac{1}{n^3}\right) \right)^p.
\end{eqnarray*} Therefore, \begin{eqnarray*}
\frac{f(n)}{(2n)^{p}} &=& \frac{1}{2^p  n^{\frac{p}{2}}} \left( 1 - \frac{1}{2\sqrt{n}} + O\left(\frac{1}{n}\right) \right)^p + \left( 1 - \frac{5}{4n} + \frac{1}{16n^2} + O\left(\frac{1}{n^4}\right) \right)^p + O\left(\frac{1}{n^p}\right)
\\&&- \frac{1}{2^p n^{\frac{p}{2}}} \left( 1 - \frac{1}{\sqrt{n}} + O\left(\frac{1}{n}\right) \right)^p - \left( 1 - \frac{5}{4n} + \frac{3}{16n^2} + O\left(\frac{1}{n^4}\right) \right)^p
\\&=& \frac{1}{2^p  n^{\frac{p}{2}}} \left[ 1 + p \left( - \frac{1}{2\sqrt{n}} \right) + O\left(\frac{1}{n}\right) \right]
+ 1 + p \left( - \frac{5}{4n} + \frac{1}{16n^2} + O\left(\frac{1}{n^4}\right) \right)
\\&&+ \frac{p(p-1)}{2} \left( -\frac{5}{4n} +\frac{1}{16n^2} + O\left(\frac{1}{n^4}\right)\right)^2+O\left(\frac{1}{n^3}\right)+O\left(\frac{1}{n^p}\right)
\\&&- \frac{1}{2^p n^{\frac{p}{2}}} \left[ 1 + p \left( - \frac{1}{\sqrt{n}} + O\left(\frac{1}{n}\right) \right) \right] -  1 - p \left( - \frac{5}{4n} + \frac{3}{16n^2} + O\left(\frac{1}{n^4}\right) \right) \\&&- \frac{p(p-1)}{2} \left( -\frac{5}{4n}+ \frac{3}{16n^2}+ O\left(\frac{1}{n^4}\right)\right)^2 + O\left(\frac{1}{n^3}\right)
\\&=& \frac{p}{2^{p+1} n^{\frac{p+1}{2}}} - \frac{p}{8n^2} + O\left(\frac{1}{n^{\frac{p}{2}+1}}\right)+O\left(\frac{1}{n^p}\right),
\end{eqnarray*} for \( 1 < p < 3 \). Notice that \( \frac{p+1}{2} < 2 \), \( \frac{p+1}{2} < \frac{p}{2}+1 \), and \( \frac{p+1}{2} < p \), there exists a positive integer \( n_0 \) such that \( \frac{f(n)}{(2n)^p} > 0 \) holds for all \( n \geq n_0 \). In other words, for sufficiently large \( n \), we have \( \mathcal{E}_p^{+}\left( \overline{S_{n,n}} \right) > \mathcal{E}_p^{+}\left( (\overline{S_{n,n}})^+ \right) \). This completes the proof. \end{proof}

\begin{remark}
By setting \( p = 2 \) in Theorem \ref{t5}, we can directly obtain a family of counterexamples to Conjecture \ref{conj2}.
\end{remark}

\section{Concluding Remarks and Open Problems}\label{sec4}
In the previous sections, we introduced the concepts of positive and negative $p$-energies and derived some of their properties. We believe this extension paves the way for numerous unexplored areas that merit further investigation.

In the proof of Theorem~\ref{t5}, we selected the pair of graphs \( \overline{S_{n,n}} \) and \( (\overline{S_{n,n}})^+ \). It is worth noting that this type of construction in fact generalizes the example presented in~\cite[Example~2.6]{Abiad2023}. In that example, their graph \( G - e \) is precisely \( \overline{S_{5,4}} \), and \( G \) is \( (\overline{S_{5,4}})^+ \). Through similar computations, one can verify that the pair \( \overline{S_{n+1,n}} \) and \( (\overline{S_{n+1,n}})^+ \) can also be used to prove Theorem~\ref{t5}. However, since computing the eigenvalues of graphs with odd order is generally less convenient than for those of even order, we chose to focus on extremal graphs of even order in Theorem~\ref{t5}.

Unfortunately, whether we use the odd-order graphs \( \overline{S_{n+1,n}} \) or the even-order graphs \( \overline{S_{n,n}} \), neither construction is able to refute the monotonicity property for \( p \geq 3 \) in Question~\ref{q1}. It is possible that these cases exhibit monotonicity with respect to edge addition. Therefore, we reformulate Guo's conjecture as the following:

\begin{prob}\label{conjsr}
Let $p \geq 3$. If $G$ is a graph and $uv \notin E(G)$, does it always hold that
\[
\mathcal{E}_p^+(G + uv) \geq \mathcal{E}_p^+(G)\,?
\] 
\end{prob}

Define $T_p(n):= \min \left\{ \mathcal{E}_p(G)|\, G \text{ is an } n \text{-vertex tree} \right\}$, and it is clear that $T_2(n) = 2(n-1)$ since all trees have $n-1$ edges. Therefore, assume Conjecture \ref{conjsr} is true, then for any connected graph \( G \), we have \( \mathcal{E}_p^+(G) \geq \mathcal{E}_p^+(T_G) \), where \( T_G \) is a spanning tree of \( G \). Combining with the well-known fact that if $\lambda$ is an eigenvalue of a tree $T$, then $-\lambda$ is also an eigenvalue \cite[Lemma 1]{Cao1995}, we know that \(
\mathcal{E}_p^+(G) \geq \frac{1}{2}T_p(n)\) holds for any connected graph $G$ when $p \geq 3$.

At this point, one may naturally ask whether the exact value of \( T_p(n) \) can be determined. This remains an open conjecture, with only partial results known so far. Specifically, let \( P_n \) denote the path graph on \( n \) vertices, and \( S_n \) the star graph consisting of one central vertex connected to \( n - 1 \) leaves. It is well known that, among all trees with \( n \) vertices, the star graph \( S_n \) and the path graph \( P_n \) attain the minimum and maximum graph energy, respectively. However, for the \( p \)-energy, the extremal behavior depends on the value of \( p \). In particular, different patterns emerge in the ranges \( 1 \leq p \leq 2 \) and \( p > 2 \). This phenomenon was first observed by S.~Wagner (as cited in a private communication in~\cite{Li2012}) and was later reiterated by Nikiforov~\cite{Nikiforov2016}, who formulated the following conjecture.

\begin{conj}[\cite{Nikiforov2016}]\label{conj418}
Let $T = (V, E)$ be a tree with $|V| = n$ and $|E| = m$.
\begin{itemize}
  \item[(i)] If $1 \leq p \leq 2$, then $\mathcal{E}_p(S_n) \leq \mathcal{E}_p(T) \leq \mathcal{E}_p(P_n)$.
  \item[(ii)] If $p > 2$, then $\mathcal{E}_p(P_n) \leq \mathcal{E}_p(T) \leq \mathcal{E}_p(S_n)$.
\end{itemize}
\end{conj}

The first result concerns the case \( 1 \leq p \leq 2 \). In~\cite[Theorem 2]{Arizmendi23}, Arizmendi and Guerrero proved the following theorem, which confirms Conjecture~\ref{conj418}~(i).

\begin{thm}[\cite{Arizmendi23}]
Let \( S_n \) be the star graph, \( P_n \) the path graph, and \( T_n \) a tree on \( n \) vertices. Then, for \( 0 < p \leq 2 \),
\[
\mathcal{E}_p(S_n) \leq \mathcal{E}_p(T_n) \leq \mathcal{E}_p(P_n).
\]
\end{thm}

The second result concerns the case \( p > 2 \). Arizmendi and Arizmendi~\cite[Proposition~4.7~(ii)]{Arizmendi2023} confirmed the right-hand side of Conjecture~\ref{conj418}~(ii):

\begin{thm}[\cite{Arizmendi2023}]
Let \( S_n \) be the star graph, and let \( T_n \) be a tree on \( n \) vertices. Then, for \( p > 2 \),
\[
\mathcal{E}_p(T_n) \leq \mathcal{E}_p(S_n).
\]
\end{thm}

The final result concerns the case where \( p > 2 \) is a positive even integer. In~\cite{Csikvari2010}, Csikv\'{a}ri proved that among all connected graphs of a given order, the path graph minimizes the number of closed walks of any given length. This result was later reformulated by Nikiforov~\cite[Proposition 4.49]{Nikiforov2016} as follows:

\begin{thm}[\cite{Nikiforov2016}]
Let \( G \) be a connected graph of order \( n \). Then, for every integer \( k \geq 2 \), we have
\[
\mathcal{E}_{2k}(G) \geq \mathcal{E}_{2k}(P_n).
\]
\end{thm}

Motivated by the evidence presented above, we propose a conjecture at the end of this paper that generalizes Conjecture~\ref{conj:extension-Hong}, specifically concerning the positive $p$-energy \( \mathcal{E}_p^{+}(G) \), stated as follows:

\begin{conj}\label{conjpositive}
Let \( p \geq 2 \), and let \( G \) be a connected graph with \( n \) vertices. Then
\[
\mathcal{E}_p^+(G) \geq \mathcal{E}_p^+(P_n) = \sum_{k=1}^{\left\lfloor \frac{n+1}{2} \right\rfloor} 2^p \cos^p\left( \frac{k \pi}{n+1} \right).
\]
\end{conj}

\begin{remark}
We have examined all connected graphs with up to 10 vertices and tested various values of \( p \). No counterexamples to Conjecture~\ref{conjpositive} have been found.
\end{remark}

\begin{remark}
If the answers to Problem~\ref{conjsr} and Conjecture~\ref{conj418} are both affirmative, they can be used to establish the case \( p \geq 3 \) of Conjecture~\ref{conjpositive}. Indeed, this would imply \(\mathcal{E}_p^+(G) \geq \mathcal{E}_p^+(T_G) \geq \mathcal{E}_p^+(P_n)\), where \( T_G \) denotes a spanning tree of \( G \). Therefore, for large values of \( p \), edge-addition monotonicity could serve as a plausible approach to proving Conjecture~\ref{conjpositive}.
\end{remark}

\section*{Acknowledgments}
The authors are grateful to Prof.~Aida Abiad, Dr.~Clive Elphick, and Prof.~Krystal Guo for their valuable comments and suggestions, which significantly improved an earlier version of this paper. The research of the third author is supported by National Key Research and Development Program of China 2023YFA1010203 and National Natural Science Foundation of China (Grant No.\,12371357).


\begin{thebibliography}{99}

\bibitem{Abiad2023}
A. Abiad, L. Lima, D. N. Desai, K. Guo, L. Hogben, and J. Madrid,
Positive and negative square energies of graphs, \textit{Electron. J. Linear Algebra}, \textbf{39} (2023): 307--326.

\bibitem{Akbari2020}
S. Akbari, M. Einollahzadeh, M.M. Karkhaneei, and M.A. Nematollahi, Proof of a conjecture on the Seidel energy of graphs, \textit{Eur. J. Comb.}, \textbf{86} (2020): 103078.

\bibitem{Akbari2024}
S. Akbari, H. Kumar, B. Mohar, and S. Pragada, A linear lower bound for the square energy of graphs, \textit{arXiv preprint}, arXiv:2409.18220, 2024.

\bibitem{Akbari2025}
S. Akbari, H. Kumar, B. Mohar, and S. Pragada, Vertex partitioning and \( p \)-energy of graphs, \textit{arXiv preprint}, arXiv:2503.16882, 2025.

\bibitem{Ando2015}
T. Ando and M. Lin, Proof of a conjectured lower bound on the chromatic number of a graph, \textit{Linear Algebra Appl.}, \textbf{485} (2015): 480--484.

\bibitem{Arizmendi2023}
G. Arizmendi, O. Arizmendi, The graph energy game, \textit{Discrete Appl. Math.}, \textbf{330} (2023): 128–140.

\bibitem{Arizmendi23}
O. Arizmendi and J. Guerrero, On the $p$-Schatten energy of bipartite graphs, \textit{Acta Math. Hungar.}, \textbf{169(2)} (2023): 503--509.




\bibitem{Cao1995}
D. Cao and H. Yuan, The distribution of eigenvalues in graphs, \textit{Linear Algebra Appl.}, \textbf{216} (1995): 211--224.


\bibitem{Cioaba}
S.~M. Cioaba, K. Guo, and N. Srivastava, Spectral graph and hypergraph theory: connections and applications, 2021.


\bibitem{Csikvari2010}
P. Csikv\'{a}ri, On a poset of trees, \textit{Combinatorica}, \textbf{30} (2010): 125--137.



\bibitem{Cvetkovic2010}
D. Cvetkovi\'{c}, P. Rowlinson, and S.~K. Simi\'{c}, \textit{An Introduction to the Theory of Graph Spectra}, Cambridge University Press, 2010.



\bibitem{EFGW16}
C. Elphick, M. Farber, F. Goldberg, and P. Wocjan,
Conjectured bounds for the sum of squares of positive eigenvalues of a graph, \textit{Discrete Math.}, \textbf{339} (2016): 2215--2223.


\bibitem{ElphickLinz2024}
C. Elphick and W. Linz, Symmetry and asymmetry between positive and negative square energies of graphs, \textit{Electron. J. Linear Algebra}, \textbf{40} (2024): 418--432.

\bibitem{Elphick2017}
C. Elphick and P. Wocjan, An inertial lower bound for the chromatic number of a graph, \textit{Electron. J. Combin.}, \textbf{24(1)} (2017): \#P1.58.

\bibitem{Hall2009}
F.~J. Hall, K. Patel, and M. Stewart, Interlacing results on matrices associated with graphs, \textit{Combinatorial Mathematics and Combinatorial Computing}, \textbf{68} (2009): 113--127.


\bibitem{Horn13}
R.~A. Horn and C.~R. Johnson, \textit{Matrix Analysis}, 2nd ed., Cambridge University Press, 2013.

\bibitem{Li2012}
X. Li, Y. Shi and I. Gutman, Graph Energy, \textit{Springer}, New York, 2012.

\bibitem{LiuNing2023}
L. Liu and B. Ning, Unsolved problems in spectral graph theory, \textit{Oper. Res. Trans.}, \textbf{27(4)} (2023): 33--60.


\bibitem{Lovasz1973}
L. Lov\'{a}sz and J. Pelik\'{a}n, On the eigenvalues of trees, \textit{Period. Math. Hungar.}, \textbf{3} (1973): 175--182.



\bibitem{MarshallOlkin1979}
A.~W. Marshall and I. Olkin, \textit{Inequalities: Theory of Majorization and Its Applications}, Academic Press, 1979.




\bibitem{Nikiforov2012}
V. Nikiforov, Extremal norms of graphs and matrices, \textit{J. Math. Sci.}, \textbf{182(2)} (2012): 164--174.




\bibitem{Nikiforov2016}
V. Nikiforov, Beyond graph energy: Norms of graphs and matrices, \textit{Linear Algebra Appl.}, \textbf{506} (2016): 82-138.

\bibitem{Tang2025}
Q. Tang and C. Elphick, A spectral lower bound on the chromatic number using $p$-energy, \textit{arXiv preprint}, arXiv:2504.01295, 2025.


\bibitem{Zhang2024}
S. Zhang, Extremal values for the square energies of graphs, \textit{arXiv preprint}, arXiv:2409.15504, 2024.

\end{thebibliography}
\end{document}